**The finite precision computation cause the nonconvergence of difference scheme**

Wang Pengfei[1] , Li Jianping

State Key Laboratory of Numerical Modeling for Atmospheric Sciences and Geophysical Fluid Dynamics (LASG), Institute of Atmospheric Physics, Chinese Academy of Sciences, Beijing, 100029, China

**Abstract** The authors show that the round-off error can break the consistency which is the premise of using the difference equation to replace the original differential equations. We therefore proposed a theoretical approach to investigate this effect, and found that the difference scheme can not guarantee the convergence of the actual compute result to the analytical one. A conservation scheme experiment is applied to solve a simple linear differential equation satisfing the LAX equivalence theorem in a finite precision computer. The result of this experiment is not convergent when time step-size decreases trend to zero, which proves that even the stable scheme can't guarantee the numerical convergence in finite precision computer. Further the relative convergence concept is introduced.

## **1** Introduction

The study of stability and convergence are often connected with the LAX equivalence theorem (from here LAX theorem). This theorem can be described as: 'given a properly posed initial-value problem and a finite-difference approximation to it that satisfies the consistency condition, stability is the necessary and sufficient condition for convergence [1-2]'. The theorem is first proved with linear equations and explicit difference scheme. Richtmyer [2] applied it with implicit scheme for the linear equations. Henrici [3] found that 'the scheme is convergent if and only if it is both stable and consistent' is true even for some nonlinear cases. Rosinger[4] gives an attempt to extend LAX theorem to nonlinear case using semigroup transform. Because it is easier to obtain and investigate stability for a differential equation than to obtain the convergence, so many researchers use the stability scheme to obtain the numerical solution without the analysis of the convergence, and assume that the scheme and solution is convergent unconditionally through LAX theorem.

[1] Corresponds author Wang pengfei.
email address wpf@mail.iap.ac.cn

Lax and Richtmyer had realized that the round-off error in the computation may affect the stability and convergence problem, but for the convenience of investigation they didn't consider the round-off errors in their research (paragraph 2 of paper [1]). So we must be aware that the LAX theorem is obtained from the theoretical analysis of numerical mathematics rather than from the actual computation though the theorem is almost correct in most cases.

The study of round-off error in numerical computation can go all the way back to the time before the modern computer was invented. It was discussed by astronomers [5,6] then. The pioneering important work on the analysis of numerical error with round-off error can be found in the Neumann[7] and Turing's[8] paper soon after the first computer was invented. Mitchell discussed the round-off error difference method [9,10] and later Wilkinson[11] and Henrici[12,3] investigated the round-off error in algebraic and difference process. The more comprehensive introduction of round-off study can be found in the book by Higham[13] and the reference cited therein. Most discussions of round-off error are about how they cause the shortage of stability and convergence etc, and the behavior is still far beyond analytical analysis.

The studies in the late 20$^{th}$ century indicate that the round-off error may have effect beyond our expectation to the computation results. Li [14] et al.'s experiments showed that single- and double-precision floating point operations have important effects on the long-time numerical integration in nonlinear systems. They identified the Optimal Step-size (OS) and Maximum Effective Computation Time (MECT) using an optimal searching method. Moreover, Li[15] et al. obtained the formulas of OS and MECT through theoretical analysis. They used the improved prior bounds of discretization error to discuss the estimation of ordinary differential equations' error boundary, and obtained the relationship between OS and computation precision and the order of the method. Wang[16] et al improved their experiments by using multiple precisions to conduct further analysis for nonlinear equations.

In the following section we first investigate the nonconvergency by theoretical analysis and then present a experiment to validate the theory. Further we propose the relative convergence concept in finite precision computation.

## 2 Theory for the noconvergence

Designate $A$ as a linear operation (**the nonlinear case is analysis in appendix**) that transforms the elements $u$ into the element $Au$ by matrix-vector multiplication. The initial value problem of differential equation is as follow:

$$\begin{aligned} &\frac{d}{dt}u(t) = Au(t) \\ &u(0) = u_0 \end{aligned} \tag{1}$$

To solve the equation, the difference scheme is applied as :

$$\begin{aligned} &\frac{u^{n+1} - u^n}{\Delta t} = Au^n \\ &u^0 = u_0 \end{aligned} \tag{2}$$

Where $u^n$ denotes the numerical solution of step $n$ and $\Delta t$ is the step-size.

The **classical convergence** only considers the discreatization error, and defined by:

$$\left\| u^n - u \right\| \to 0 \tag{3}$$

When the step-size $\Delta t \to 0$

The **consistency condition** is:

$$\lim_{\Delta t \to 0} \left\| \frac{u^{n+1} - u^n}{\Delta t} - Au \right\| = 0 \tag{4}$$

While dealing with the long time integration of time-dependent differential equations, there are two types of stabilit. One is the behavior of the solution with the mesh size trend to zero within a fixed time $T$ . Another is the solution with fixed mesh size and infinite time trend. The first issue is often regarded as the 'Classical' stability from LAX. In the some research areas such as the fluid mechanics, the astronomy orbit integral, and the weather forecasting, this second stability problem often occurs.

The stability guarantees that the numerical solution does not amplify to infinite when the integration time increases. The conservation scheme is one of the ways to keep the stability in the conservation systems. Besides this scheme property the Neumann stability condition and CFL are also required.

The **stable condition** of scheme (2) is that $u^n$ is universal bounded.

$$\left\| u^n \right\| \le K \tag{5}$$

because

$$u^{n+1} = u^n + \Delta t A u^n = \left( I + \Delta t A \right) u^n = \left( I + \Delta t A \right)^n u^0 \tag{6}$$

So the stable condition can be regard as:

$$\left\| \left( I + \Delta t A \right)^n \right\| \le K \tag{7}$$

Where $K$ is a constant independent with $\Delta t$ .

But when the round-off error exists in the computation the scheme (2) should be changed to:

$$u^{n+1} = u^n + \Delta t A u^n + \varepsilon^n \tag{8}$$

Where $\varepsilon^n$ is the round-off error in one computation step[11].

Write it back to the difference format

$$\frac{u^{n+1} - u^n}{\Delta t} = A u^n + \frac{\varepsilon^n}{\Delta t} \tag{9}$$

From this formula we know that when $\Delta t \to 0$

$$\lim_{\Delta t \to 0} \left\| \frac{u^{n+1} - u^n}{\Delta t} - A u \right\| = \left\| \frac{\varepsilon^n}{\Delta t} \right\| \ne 0 \tag{10}$$

The consistence condition (4) is broken. For this reason the LAX theorem is not feasible here.

**Theorem.**

**The round-off error in the computation cause the scheme** (8) **nonconvergence.**

It is known that the unstable scheme cannot make the numerical result convergence to the analytical result, so when we want to prove the nonconvergence of actual computation we only need to deal with the stable scheme cases.

Designate $v^n$ as the analytical solution of equation (1), and the error between $u^n$ and $v^n$ is $e^n$. Then $u^n = v^n + e^n$.

Since $\left\| \frac{v^{n+1} - v^n}{\Delta t} - Av^n \right\| = O(\Delta t) = c^n$ is true for $0 < t < T$

Where $c^n \to 0$ while $\Delta t \to 0$.

$$u^n = u^{n-1} + \Delta t A u^{n-1} + \varepsilon^{n-1}$$
$$v^n + e^n = v^{n-1} + e^{n-1} + \Delta t A\left(v^{n-1} + e^{n-1}\right) + \varepsilon^{n-1}$$
$$e^n = e^{n-1} + \Delta t A\left(e^{n-1}\right) + \varepsilon^{n-1} + \Delta t \cdot c^{n-1}$$

From the error iterative formula

$$e^n = e^{n-1} + \Delta t A\left(e^{n-1}\right) + \varepsilon^{n-1} + \Delta t \cdot c^{n-1}$$
$$= (I + \Delta t A)e^{n-1} + \varepsilon^{n-1} + \Delta t \cdot c^{n-1}$$
$$= (I + \Delta t A)^2 e^{n-2} + (I + \Delta t A)\varepsilon^{n-2} + \varepsilon^{n-1} + (I + \Delta t A)\Delta t \cdot c^{n-2} + \Delta t \cdot c^{n-1}$$
$$= \dots$$
$$= (I + \Delta t A)^n e^0 + (I + \Delta t A)^{n-1} \varepsilon^0 + \dots + (I + \Delta t A)\varepsilon^{n-2} + \varepsilon^{n-1}$$
$$+ \Delta t \cdot \left((I + \Delta t A)^{n-1} c^0 + \dots + (I + \Delta t A)c^{n-2} + c^{n-1}\right)$$

Thus because the scheme is stable, the item $\Delta t \cdot \left((I + \Delta t A)^{n-1} c^0 + \dots + (I + \Delta t A)c^{n-2} + c^{n-1}\right) \to 0$

thus

$$e^n = (I + \Delta t A)^n e^0 + (I + \Delta t A)^{n-1} \varepsilon^0 + \dots + (I + \Delta t A)\varepsilon^{n-2} + \varepsilon^{n-1}$$

Because $\varepsilon^n$ is independent with $\Delta t$, it is generally not uniform convert to 0 while $\Delta t \to 0$. Thus the theorem is finished.

This result can be compared to the result of Bruno[17], in whose paper the perturbation term is $\varepsilon^n = \Delta t \cdot s^n = O(\Delta t)$. The perturbation term in our study is $\varepsilon^n = O(1)$. This difference causes the different convergence behavior to the same difference scheme.

For the implicit scheme, it can be transfered to an explicit scheme after solving the linear algebraic equations. Thus it also conforms to the nonconvergence behavior in finite precision computation.

The above discussion is fit for the partial differential $\frac{\partial}{\partial t}u = Au$ where $A$ is linear operator (the discussion of nonlinear case is in appendix).

## 3 Nonconvergences: the experiments validation

To test if the round-off error can really cause the nonconvergence in actual computation, we apply Euler midpoint scheme with a simple equation to compute the actual solution. The 3rd order conservation Runge-kutta scheme is also applied which is detailed in the Wang[18].

### 3.1 The equation and difference scheme:

We can obtain conservation scheme to solve the equation:

$$\begin{cases} \frac{\mathrm{d}x}{\mathrm{d}t} = -ay \\ \frac{\mathrm{d}y}{\mathrm{d}t} = bx \end{cases} \tag{11}$$

As we know the analytical solution is (in our study $a = 0.1, b = 0.2$):

$$\begin{cases} x = \cos\left(\sqrt{ab}t\right) \\ y = \sqrt{\frac{b}{a}}\sin\left(\sqrt{ab}t\right) \end{cases} \tag{12}$$

This equation has been used to analyze numerical error before[19], but the analysis did not focus on the convergence discussion.

The Euler mid-point scheme is an implicit scheme:

$$\frac{F^{n+1} - F^{n}}{\Delta t} + A\left(\frac{F^{n+1} + F^{n}}{2}\right) = 0 \tag{13}$$

Convert to equations:

$$\begin{cases} \dfrac{x^{n+1} - x^{n}}{\Delta t} + a\left(\dfrac{y^{n+1} + y^{n}}{2}\right) = 0 \\ \dfrac{y^{n+1} - y^{n}}{\Delta t} - b\left(\dfrac{x^{n+1} + x^{n}}{2}\right) = 0 \end{cases} \tag{14}$$

The solution is:

$$\begin{cases} x^{n+1} = \dfrac{x^{n}\left(1 - \dfrac{a\Delta t}{2}\dfrac{b\Delta t}{2}\right) - a\Delta t y^{n}}{1 + \dfrac{a\Delta t}{2}\dfrac{b\Delta t}{2}} \\ y^{n+1} = \dfrac{y^{n}\left(1 - \dfrac{a\Delta t}{2}\dfrac{b\Delta t}{2}\right) + b\Delta t x^{n}}{1 + \dfrac{a\Delta t}{2}\dfrac{b\Delta t}{2}} \end{cases} \tag{15}$$

It is easy to validate the inner product conservation since

$$\left(AF^{n+1}, F^{n+1}\right) = bx^{n+1}x^{n+1} + ay^{n+1}y^{n+1} =$$

$$b\left(\frac{x^{n}\left(1 - \dfrac{a\Delta t}{2}\dfrac{b\Delta t}{2}\right) - a\Delta t y^{n}}{1 + \dfrac{a\Delta t}{2}\dfrac{b\Delta t}{2}}\right)^{2} + a\left(\frac{y^{n}\left(1 - \dfrac{a\Delta t}{2}\dfrac{b\Delta t}{2}\right) + b\Delta t x^{n}}{1 + \dfrac{a\Delta t}{2}\dfrac{b\Delta t}{2}}\right)^{2}$$

$$= bx^{n}x^{n} + ay^{n}y^{n}$$

Write it in matrix format:

$$\begin{pmatrix} x^{n+1} \\ y^{n+1} \end{pmatrix} = \begin{pmatrix} \dfrac{\left(1-\dfrac{a\Delta t}{2}\dfrac{b\Delta t}{2}\right)}{1+\dfrac{a\Delta t}{2}\dfrac{b\Delta t}{2}} & -\dfrac{a\Delta t}{1+\dfrac{a\Delta t}{2}\dfrac{b\Delta t}{2}} \\ \dfrac{b\Delta t}{1+\dfrac{a\Delta t}{2}\dfrac{b\Delta t}{2}} & \dfrac{\left(1-\dfrac{a\Delta t}{2}\dfrac{b\Delta t}{2}\right)}{1+\dfrac{a\Delta t}{2}\dfrac{b\Delta t}{2}} \end{pmatrix} \begin{pmatrix} x^{n} \\ y^{n} \end{pmatrix} = A \begin{pmatrix} x^{n} \\ y^{n} \end{pmatrix}$$

Where $A$ is matrix, $A = \begin{pmatrix} \dfrac{\left(1-\dfrac{a\Delta t}{2}\dfrac{b\Delta t}{2}\right)}{1+\dfrac{a\Delta t}{2}\dfrac{b\Delta t}{2}} & -\dfrac{a\Delta t}{1+\dfrac{a\Delta t}{2}\dfrac{b\Delta t}{2}} \\ \dfrac{b\Delta t}{1+\dfrac{a\Delta t}{2}\dfrac{b\Delta t}{2}} & \dfrac{\left(1-\dfrac{a\Delta t}{2}\dfrac{b\Delta t}{2}\right)}{1+\dfrac{a\Delta t}{2}\dfrac{b\Delta t}{2}} \end{pmatrix}$.

It is known that $\begin{pmatrix} x^{n+1} \\ y^{n+1} \end{pmatrix} = A \begin{pmatrix} x^{n} \\ y^{n} \end{pmatrix} = A^{n} \begin{pmatrix} x^{0} \\ y^{0} \end{pmatrix}$ is stable when the spectral radius of $A$, $\rho(A) < 1$ in the exact arithmetic. In this case when $\Delta t \to 0$, the condition $\rho(A) < 1$ is satisfied.

F. Chatelin [20]investigated the convergence of linear iterative scheme:

$$x^{k+1} = Ax^{k} + c \qquad (16)$$

, and found that when $\rho(A) < 1$ due to the finite precision the computed norm have three possible results: (a) for a small nonnormality $\|A^{k}\|$ behaves like it in the exact arithmetic, (b) for a moderate nonnormality $\|A^{k}\|$ oscillates, (c) for a large nonnormality $\|A^{k}\|$ diverges.

In our case the perturbation $c = \varepsilon^{n}$, the computed convergence result behavior can be obtained from experiments.

### 3.2 Error separation formula

To evaluate the error between computed solution and the theoretical solution we should bring in the error formula. The solution is sited in an ellipse as shown in figure

1.

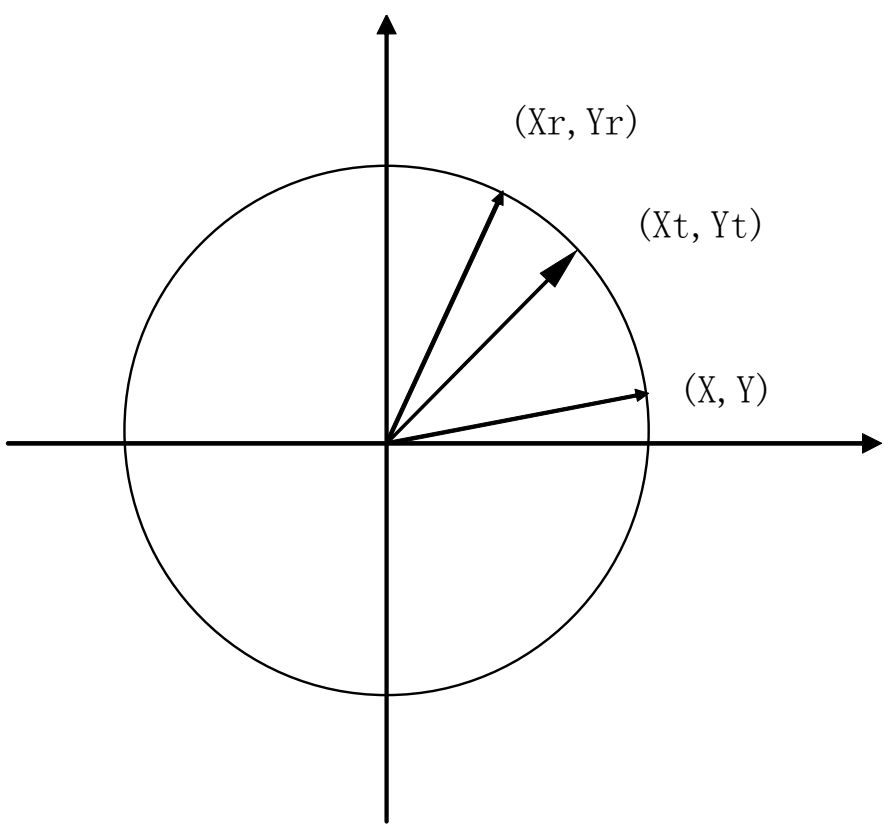

Fig. 1. The demonstration of the theoretical solution $(X,Y)$, the reference solution $(X_t,Y_t)$, and the actual solution $(X_r,Y_r)$.

The notation $(X,Y)$ is the theoretical solution, the $(X_t,Y_t)$ is the solution which only truncation error in it (usually we can get a reference solution close to the solution of exact arithmetic), and the $(X_r,Y_r)$ is the computed solution. Three types of error are defined here for the variable $X$.

The total error:

$$E_x = X_r - X \,, \tag{17}$$

The truncation error

$$E_{xt} = X_t - X \tag{18}$$

And the round-off error

$$E_{xr} = X_r - X_t \tag{19}$$

these equations (17-19) can be established as

$$E_x = E_{xr} + E_{xt} \quad . \tag{20}$$

The error of another variable of the equation $Y$ can be written as the above format too. When we need to evaluate the integrate error, we can define the norm of error.

The norm of total error is

$$E = \sqrt{\left(X_r - X\right)^2 + \left(Y_r - Y\right)^2} \tag{21}$$

the norm of truncation error

$$E_t = \sqrt{\left(X_t - X\right)^2 + \left(Y_t - Y\right)^2} \tag{22}$$

the norm of round-off error

$$E_r = \sqrt{\left(X_r - X_t\right)^2 + \left(Y_r - Y_t\right)^2} \tag{23}$$

**3.3 experiments result**

The computer arch we run the program is IBM-P690 and PC machine. A corresponding quadruple precision result is act as a reference solution which is close to exact arithmetic numerical solution. Some detail discussion of reference solution can be found at Wang[16]. In out experiments the round-off error result is computed by single-precision and the reference solution is computed bye quadruple precision.

The first experiment is integrated to $T = 10000$, and the step-size varies from 0.1 to 0.0000001.

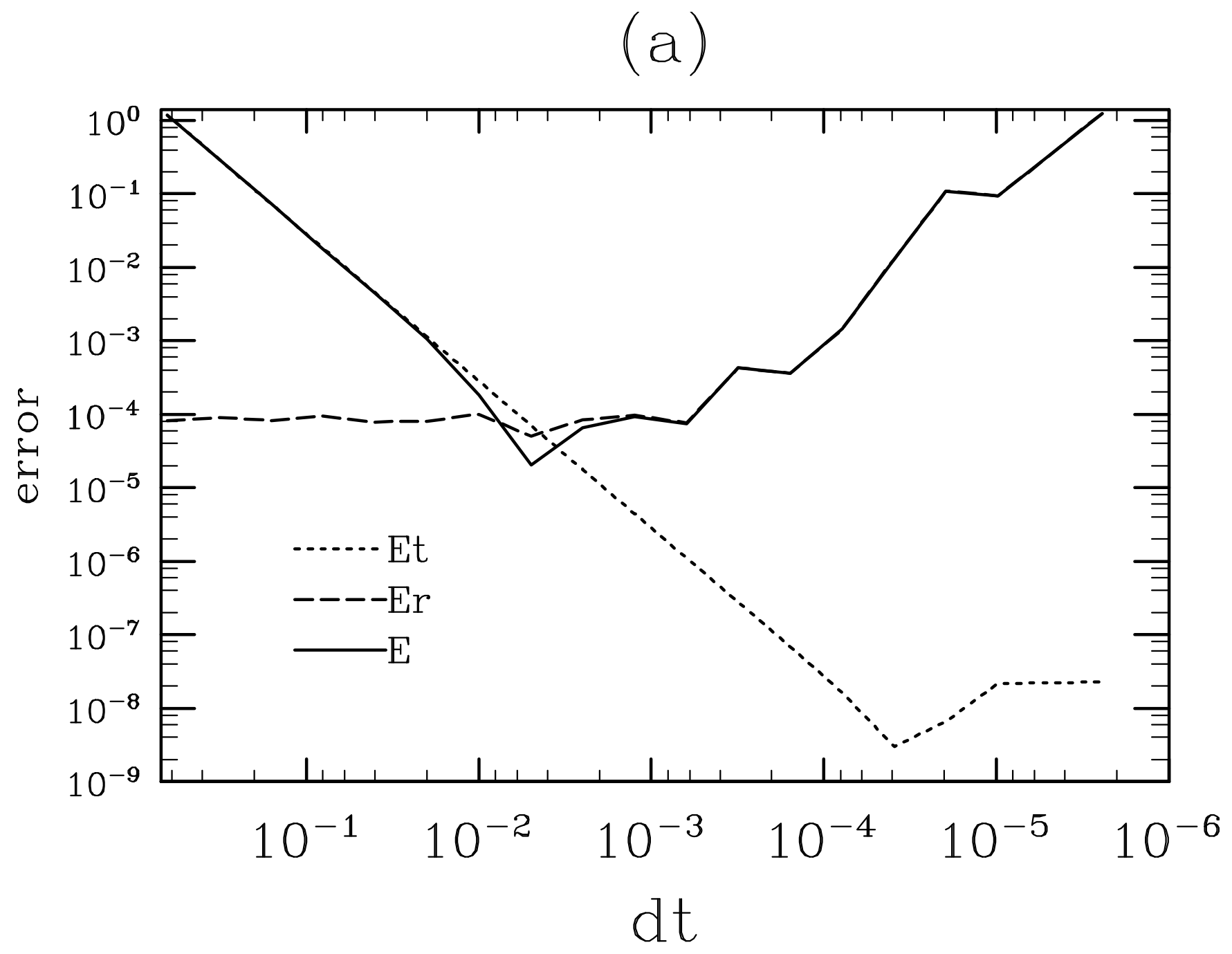

Fig. 2. The numerical solution error of time $T$=10000.0 as the step-size $\Delta t$ changes,.

As shown in figure 2, the error of reference solution $E_t$ decreases while the step-size decreases as the expected. The round-off error $E_r$ is at a small value when the $\Delta t > 10^{-4}$, but it increase when the $\Delta t$ keeps to decrease to about $10^{-6}$. The total error $E$ decreases when $\Delta t > 10^{-3}$ but begins to increase when $\Delta t < 10^{-4}$, and get a minimize value when $\Delta t \approx 10^{-3}$ .From the experiment we validate the nonconvergence analysis.

The second experiment is to do long time integral with constant $\Delta t$ such as $\Delta t = 0.00001$, $t$ varies from 0 to 100000, to investigate the error variety of $E_r$ and $E_t$. Because when $\Delta t > 10^{-3}$ the discretezation error is the main error source in total error, and it is widely discussed in many numerical analysis books, so we will focus on the error behavior of $\Delta t < 10^{-4}$ .

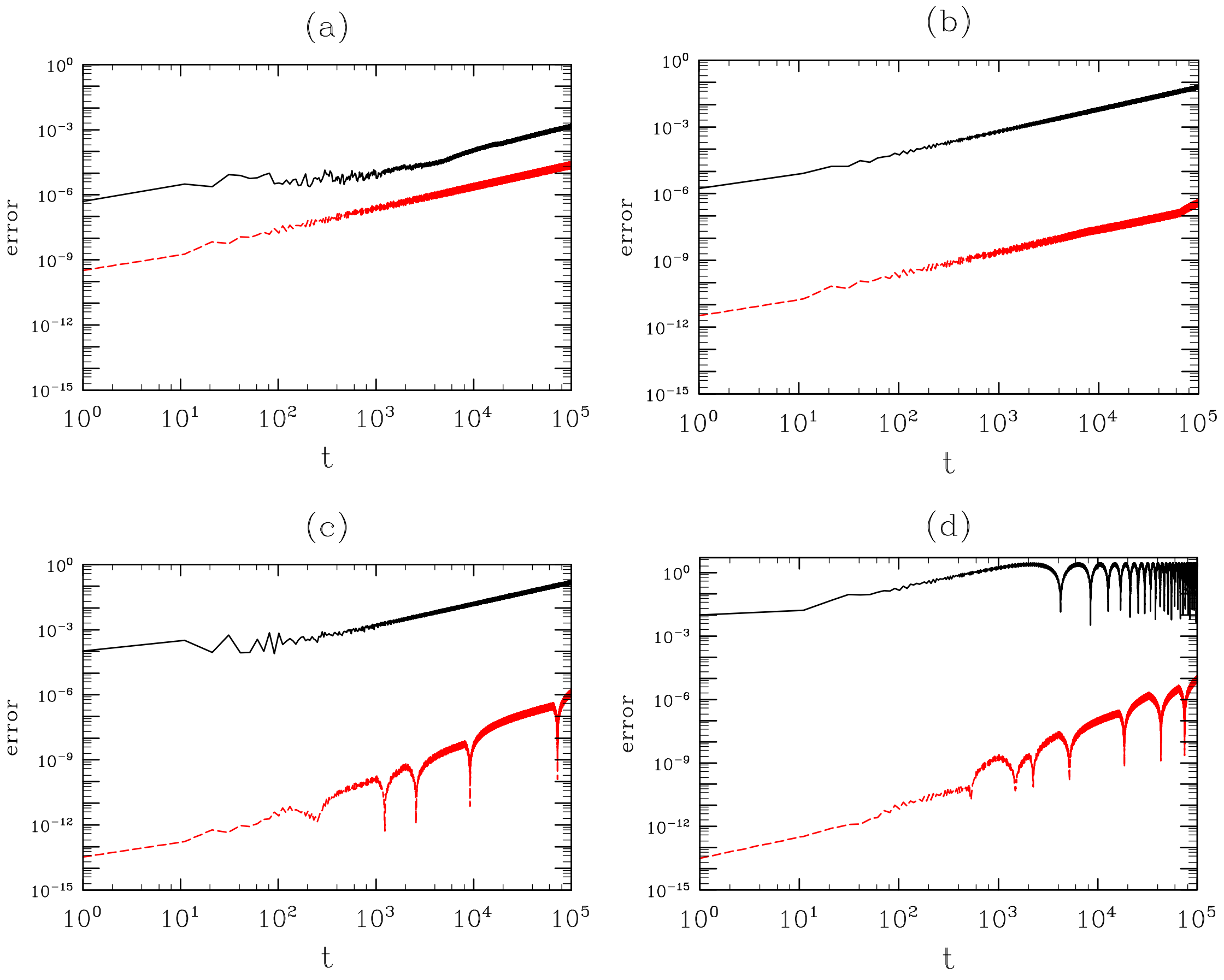

Fig. 3. The numerical solution error $E_r$ (black) and $E_t$ (red) versus time (a) $\Delta t = 10^{-3}$,

(b) $\Delta t = 10^{-4}$, (c) $\Delta t = 10^{-5}$, (d) $\Delta t = 10^{-6}$.

Figure 3 shows the behavior round-off error $E_r$ and the truncation error $E_t$ versus time $t$. When $\Delta t = 10^{-3}$ the round-off is larger than truncation error, but there is not so much difference in magnitude, and the error increases slowly while the time increases. At time $t = 10^5$ the error is about $10^{-3}$. From the figures b, c, and d we can find that when $\Delta t$ decreases the round-off error is becoming larger, on the contrary the truncation error is decrease. Especially when the $\Delta t = 10^{-6}$ the round-off error became saturation at the time $t = 10^4$. We can image that as while as the time $t$ increase, the error in (a, b, c) should be saturation too.

The experiment proved that the actual results are not uniform convergence to the analytical one while $\Delta t \to 0$ , but it seemed that the results are still bounded and stable.

## 4 The concept of relative convergent

The theoretical analysis and the numerical experiments indicate that the classical absolute convergence is not true in the real computation environments. So what convergence concept can be used to replace the absolute convergence is the question.

The relative convergence thus come into view. The focus of relative convergence is not the numerical solution convergence to the analytical solution as $\Delta t \to 0$, but to keep the total error in an admissible bound. When the total error reached the error bound we define it as the effective computation time (ECT). And the numerical solution beyond ECT is not credible from convergence view.

The relative convergent depends on the ODEs, the order of scheme and the float-point precision. In the cases where $\Delta t$ is not very small and the accumulation of round-off error is not larger than the truncation error, we can use classical numerical analysis knowledge to analyze the total error. On the contrary when the round-off error became the primary, we must use statistical analysis of round-off error or the reference solution method to determine the total error.

## 5 Conclusion

We present that the finite precision computation in difference scheme causes the absolute convergence loss. The theoretical evidence is proposed and the numerical experiments are implemented to validate the evidence.

No absolute convergence can remain in the real computational environment. We then suggest using the relative convergence to replace the absolute convergence in the real computational environment.

## Appendix:

### A.1 The nonconvergence for nonlinear case

Let $v^n$ be the analytical solution of equation (1), and the error between $u^n$ and $v^n$ is $e^n$ ,and $u^n = v^n + e^n$ .

Follow the analysis of equation (8)

$$u^n = u^{n-1} + \Delta t A u^{n-1} + \varepsilon^{n-1} \Rightarrow v^n + e^n = v^{n-1} + e^{n-1} + \Delta t A\left(v^{n-1} + e^{n-1}\right) + \varepsilon^{n-1} \tag{24}$$

Since $\left\| \frac{v^{n+1} - v^n}{\Delta t} - Av^n \right\| = O\left(\Delta t\right) = c^n$ is true for $0 < t < T$

and $c^n \to 0$ while $\Delta t \to 0$ .

We can get

$$e^n = e^{n-1} + \Delta t A\left(v^{n-1} + e^{n-1}\right) - \Delta t A\left(v^{n-1}\right) + \varepsilon^{n-1} + \Delta t \cdot c^{n-1} \tag{25}$$

From Eq.**(**24**).**

$A$ is a nonlinear operator**,** generally $A\left(v+e\right) \neq A\left(v\right) + A\left(e\right)$ **.**The Lipschitz condition for the operation $\left\| A\left(v+e\right) - A\left(v\right) \right\| \le Le$ is required when expand $A\left(v^{n-1} + e^{n-1}\right)$ to Taylor series：

$$A\left(v^{n-1} + e^{n-1}\right) = Av^{n-1} + A^{'}\left(v^{n-1}\right)e^{n-1} + A^{''}\left(v^{n-1}\right)\frac{\left(e^{n-1}\right)^2}{2} + A^{'''}\left(v^{n-1}\right)\frac{\left(e^{n-1}\right)^3}{3!} + \cdots \tag{26}$$

Define a $m \times m$ order matrix operator

$$M^{n-1} = A^{'}v^{n-1} + A^{''}v^{n-1} \cdot \frac{\left(e^{n-1}\right)^1}{2} + A^{'''}v^{n-1} \cdot \frac{\left(e^{n-1}\right)^2}{3!} + \cdots \tag{27}$$

Thus

$$A\left(v^{n-1} + e^{n-1}\right) = Av^{n-1} + M^{n-1}e^{n-1} \tag{28}$$

$$A\left(v^{n-1} + e^{n-1}\right) - Av^{n-1} = M^{n-1}e^{n-1} \tag{29}$$

$\left\| A\left(v+e\right) - A\left(v\right) \right\| \le Le \Rightarrow \left\| M^{n-1}e^{n-1} \right\| \le Le$ , So $M^{n-1}$ is bounded matrix operation.

From Eq.25 and Eq.29 we get

$$e^n = \left(I + \Delta t M^{n-1}\right)e^{n-1} + \varepsilon^{n-1} + \Delta t \cdot c^{n-1} \tag{30}$$

Where $I$ is the unit operator for $m \times m$ order matrix operator.

Thus

$$\begin{aligned}
e^n &= \left(I + \Delta t M^{n-1}\right)e^{n-1} + \varepsilon^{n-1} + \Delta t \cdot c^{n-1} \\
&= \left(I + \Delta t M^{n-1}\right)\left(I + \Delta t M^{n-2}\right)e^{n-2} + \left(I + \Delta t M^{n-1}\right)\varepsilon^{n-2} + \left(I + \Delta t M^{n-1}\right)\Delta t \cdot c^{n-2} + \varepsilon^{n-1} + \Delta t \cdot c^{n-1} \\
&= \left(I + \Delta t M^{n-1}\right)\left(I + \Delta t M^{n-2}\right)e^{n-2} + \left[\left(I + \Delta t M^{n-1}\right)\varepsilon^{n-2} + \varepsilon^{n-1}\right] + \left[\left(I + \Delta t M^{n-1}\right)\Delta t \cdot c^{n-2} + \Delta t \cdot c^{n-1}\right] \\
&= \ldots \\
&= \prod_{j=0}^{n-1}\left(I + \Delta t M^{n-1-j}\right)e^0 + \sum_{k=0}^{n-1}\prod_{j=0}^{k}\left(I + \Delta t M^{k-j}\right)\Delta t \cdot c^k + \sum_{k=0}^{n-1}\prod_{j=0}^{k}\left(I + \Delta t M^{k-j}\right)\cdot \varepsilon^k
\end{aligned}$$

the item $\sum_{k=0}^{n-1}\prod_{j=0}^{k}\left(I + \Delta t M^{k-j}\right)\Delta t \cdot c^k \to 0$, and $e^0 = 0$.

Thus $e^n$ is mainly affected by the item

$$\sum_{k=0}^{n-1}\prod_{j=0}^{k}\left(I + \Delta t M^{k-j}\right)\cdot \varepsilon^k \,, \tag{31}$$

but this item does not uniformly trend to zero.

Thus $e^n \neq 0$ while $\Delta t \to 0$, and nonconvergence is proved.

**A.2 The Taylor series of $A\left(v^{n-1} + e^{n-1}\right)$**

The Taylor series of $A(v+e)$ is

$$A(v+e) = A(v) + A^{'}(v)e + A^{''}(v)\frac{(e)^2}{2} + A^{'''}(v)\frac{(e)^3}{3!} + \cdots. \tag{32}$$

For a $m$ variable system, the $v$ and $e$ are $m \times 1$ matrix, write as

$$A=\begin{pmatrix}A_1\\ \cdots\\ A_m\end{pmatrix}, v=\begin{pmatrix}v_1\\ \cdots\\ v_m\end{pmatrix}, e=\begin{pmatrix}e_1\\ \cdots\\ e_m\end{pmatrix}, v+e=\begin{pmatrix}v_1+e_1\\ \cdots\\ v_m+e_m\end{pmatrix}$$

the $A_k$ are operators. $A_k(v)=f_k(v_1,\cdots,v_m)$ are multi variable functions.

The Taylor series of a multi variable function is known as:

$$f(v_1+e_1,v_2+e_2,\cdots,v_m+e_m)=\sum_{j=0}^{\infty}\left\{\frac{1}{j!}\left[\sum_{k=1}^{m}e_k\frac{\partial}{\partial v_k}\right]^j f(v_1,v_2,\cdots,v_m)\right\}$$

So when we do Taylor expansion to matrix $A(v+e)$, we can expand each matrix element to Taylor series and then combine them together again.

$$A(v+e)=\begin{pmatrix}A_1\\ \cdots\\ A_m\end{pmatrix}(v+e)=\begin{pmatrix}A_1(v+e)\\ \cdots\\ A_m(v+e)\end{pmatrix}$$

For the function $A_1(v+e)$

$$A_1(v+e)=A_1v+A_1^{'}(v)e+A_1^{''}(v)\frac{(e)^2}{2}+A_1^{'''}(v)\frac{(e)^3}{3!}+\cdots$$

Where $A_1^{'}=\left(\frac{\partial A_1}{\partial v_1},\cdots,\frac{\partial A_1}{\partial v_m}\right)$ is an operator.

i.e.

$$A_1^{'}(v)e=\left(\frac{\partial A_1}{\partial v_1},\cdots,\frac{\partial A_1}{\partial v_m}\right)\cdot\begin{pmatrix}e_1\\ \cdots\\ e_m\end{pmatrix}$$

$$A_1^{''}=\left(\frac{\partial}{\partial v_1}\left(\frac{\partial A_1}{\partial v_1},\cdots,\frac{\partial A_1}{\partial v_m}\right),\frac{\partial}{\partial v_2}\left(\frac{\partial A_1}{\partial v_1},\cdots,\frac{\partial A_1}{\partial v_m}\right),\cdots,\frac{\partial}{\partial v_m}\left(\frac{\partial A_1}{\partial v_1},\cdots,\frac{\partial A_1}{\partial v_m}\right)\right)$$
$$=\left(\frac{\partial A_1^{'}}{\partial v_1},\cdots,\frac{\partial A_1^{'}}{\partial v_m}\right)$$

Through the iterative procedure, we can transform $A_1(v+e)$ to the formula

$$A_1(v+e)=A_1v+M_1\cdot\begin{pmatrix}e_1\\ \cdots\\ e_m\end{pmatrix}$$

Where $M_1=A_1^{'}v+A_1^{''}v\cdot\frac{(e)^1}{2}+A_1^{'''}v\cdot\frac{(e)^2}{3!}+\cdots$

Thus because $A=\begin{pmatrix}A_1\\ \cdots\\ A_m\end{pmatrix}$, we know

$A(v+e)=Av+M\cdot e$ is true where $M=\begin{pmatrix}M_1\\ \cdots\\ M_m\end{pmatrix}$

**A.3 The Lipschitz condition**

The formula (32) is established when $A_1^{'},A_1^{''},A_1^{'''}\cdots$ are all finite.

If the $A_1^{'},A_1^{''},A_1^{'''}\cdots$ have a singularity, we need another way to analysis.

Example case 1:

$$A_1(v)=f_1(v_1,v_2)=\sqrt{v_1+v_2}$$

$$A_1^{'}=\left(\frac{\partial A_1}{\partial v_1},\frac{\partial A_1}{\partial v_2}\right)=\left(\frac{1}{2\sqrt{v_1+v_2}},\frac{1}{2\sqrt{v_1+v_2}}\right)$$

The $v_1+v_2=0$ is a singularity point to $A_1^{'},A_1^{''},A_1^{'''}\cdots$, so we can not use Taylor series to do analysis.

But we know $A_1(v+e)=f_1(v_1+e_1,v_2+e_2)=\sqrt{(v_1+e_1)+(v_2+e_2)}=\sqrt{e_1+e_2}$

$$A_1(v+e)-A_1(v)=\sqrt{e_1+e_2}=\frac{e_1+e_2}{\sqrt{e_1+e_2}}=\frac{e_1}{\sqrt{e_1+e_2}}+\frac{e_2}{\sqrt{e_1+e_2}}=\left(\frac{1}{\sqrt{e_1+e_2}},\frac{1}{\sqrt{e_1+e_2}}\right)\cdot\begin{pmatrix}e_1\\ e_2\end{pmatrix}$$

thus the matrix $M_1=\left(\frac{1}{\sqrt{e_1+e_2}},\frac{1}{\sqrt{e_1+e_2}}\right)$ can still be obtained.

Since the radius of $M_1,\rho(M_1)\to\infty$, the convergence can not be determined for this analysis method.

Example case 2:

$$A_1(v) = (v_1 + v_2)^{\frac{3}{2}}$$

$$A_1(v+e) - A_1(v) = (e_1 + e_2)^{\frac{3}{2}} = \left( \frac{1}{2}\frac{(e_1+e_2)^{\frac{3}{2}}}{e_1}, \frac{1}{2}\frac{(e_1+e_2)^{\frac{3}{2}}}{e_2} \right) \cdot \begin{pmatrix} e_1 \\ e_2 \end{pmatrix}$$

$$M_1 = \left( \frac{1}{2}\frac{(e_1+e_2)^{\frac{3}{2}}}{e_1}, \frac{1}{2}\frac{(e_1+e_2)^{\frac{3}{2}}}{e_2} \right)$$

If $0 < \frac{e_1}{e_2} < \infty$

The $\rho(M)$ is still bounded.

Generally，for the case that $A_1^{'}, A_1^{''}, A_1^{'''} \cdots$ have singularity, the way to obtain matrix $M_1 = (\bar{M}_1, \bar{M}_2)$ which keeps $A_1(v+e) - A_1(v) = (\hat{M}_1, \hat{M}_2) \cdot \begin{pmatrix} e_1 \\ e_2 \end{pmatrix}$ is to set

$$\hat{M}_1 = \frac{1}{2}\frac{A_1(v+e) - A_1(v)}{e_1}$$

$$\hat{M}_2 = \frac{1}{2}\frac{A_1(v+e) - A_1(v)}{e_2}$$

If $e_1 = 0$, we then set

$$\hat{M}_1 = 0$$

$$\hat{M}_2 = \frac{A_1(v+e) - A_1(v)}{e_2}$$

The analysis that whether $\rho(M)$ is bounded can decide the convergence of the scheme. The condition can be write as $\|A_1(v+e) - A_1(v)\| \le Le$, and it is the Lipschitz condition for the function.

**A.4 Corollary**

**Theorem：**

**The stable and consistent scheme guarantee the convergence for nonlinear cases in exact computation when the operator $A$ satisfies Lipschitz condition.**

From above disscusion, when we let each $\varepsilon^k \equiv 0$ in formula (31) the error $e^n \to 0$ , thus the Corollary is then obtained.